\documentclass[letterpaper, 10 pt, conference]{ieeeconf}  % Comment this line out if you need a4paper

\IEEEoverridecommandlockouts                              % This command is only needed if 
\usepackage{graphics} % for pdf, bitmapped graphics files
\usepackage{epsfig} % for postscript graphics files
\usepackage{mathptmx} % assumes new font selection scheme installed
\usepackage{times} % assumes new font selection scheme installed
\usepackage{amsmath} % assumes amsmath package installed
\usepackage{amssymb}  % assumes amsmath package installed
\usepackage{bm}
\usepackage{color}

\usepackage{tabularx}
\usepackage{multirow}
\usepackage{rotating}
\usepackage{csquotes}
\usepackage{amsfonts}
\usepackage{xfrac}
\usepackage{empheq}
\usepackage{fancyhdr}
\usepackage{cases}
\usepackage{mathtools}
\usepackage{inputenc}
\usepackage{url}
\usepackage{fix-cm}
\usepackage{bm}
\usepackage{tocloft}
\usepackage{xcolor}
\usepackage{etoolbox}
\usepackage{latexsym}
\usepackage{textcomp}
\usepackage{setspace}
\usepackage{graphicx}
\DeclareGraphicsExtensions{.eps,.pdf,.jpg}
\usepackage{epstopdf}
\usepackage{lipsum}
\usepackage{ntheorem}

\usepackage{bibentry}

\usepackage{dsfont}

\usepackage{caption}
\usepackage{float}

\usepackage[]{algorithm2e}

\newtheorem{theorem}{Theorem}

\title{\LARGE \bf
Solving Optimal Control Problems for Delayed Control-Affine Systems with Quadratic Cost by Numerical Continuation
}

\author{Riccardo Bonalli$^{*}$ $^{\dagger}$, Bruno H\'{e}riss\'{e}$^{\dagger}$ and Emmanuel Tr\'{e}lat$^{*}$% <-this % stops a space
\thanks{$^*$Sorbonne Universit\'{e}s, UPMC Univ Paris 06, CNRS UMR 7598, Laboratoire Jacques-Louis Lions, F-75005, Paris, France.}%
\thanks{$\dagger$Onera - The French Aerospace Lab, F~-~91761 Palaiseau, France.}%
\thanks{\footnotesize riccardo.bonalli@onera.fr, bruno.herisse@onera.com,}%
\thanks{emmanuel.trelat@upmc.fr}%
}

\begin{document}

\maketitle
\thispagestyle{empty}
\pagestyle{empty}

%%%%%%%%%%%%%%%%%%%%%%%%%%%%%%%%%%%%%%%%%%%%%%%%%%%%%%%%%%%%%%%%%%%%%%%%%%%%%%%%
\begin{abstract}

In this paper we introduce a new method to solve fixed-delay optimal control problems which exploits numerical homotopy procedures. It is known that solving this kind of problems via indirect methods is complex and computationally demanding because their implementation is faced with two difficulties: the extremal equations are of mixed type, and besides, the shooting method has to be carefully initialized. Here, starting from the solution of the non-delayed version of the optimal control problem, the delay is introduced by numerical homotopy methods. Convergence results, which ensure the effectiveness of the whole procedure, are provided. The numerical efficiency is illustrated on an example.

\end{abstract}

%%%%%%%%%%%%%%%%%%%%%%%%%%%%%%%%%%%%%%%%%%%%%%%%%%%%%%%%%%%%%%%%%%%%%%%%%%%%%%%%
%\begin{eqnarray*}
%\begin{cases}
%\displaystyle \min \quad C_T(u_{\tau}(\cdot)) = \int_0^T u^2_{\tau}(t) \ dt  \ , \quad u_{\tau}(\cdot) \in L^2([0,T],\mathbb{R}) \medskip \\
%\dot{x}_{\tau}(t) = f_0(x_{\tau}(t)) + f_1(x_{\tau}(t-\tau)) + u_{\tau}(t) f_2(x_{\tau}(t)) \medskip \\
%x_{\tau}(t) = \phi(t) \ , \quad t \in [-M,0] \ \ , \qquad x_{\tau}(T) = \bar x \medskip \\
%\end{cases}
%\end{eqnarray*}optimal control problem ((\textbf{OCP})$_{\tau}$ from now on)
\section{Introduction}

Let $n$ be a positive integer, $M$, $T$ positive real numbers and define an initial state function $\phi \in C^0([-M,0],\mathbb{R}^n)$. For every $\tau \in (0,M]$, consider the following delayed single-input control-affine system in $\mathbb{R}^n$
\begin{equation} \label{dynDelay}
\begin{cases}
\dot{x}_{\tau}(t) = f_0(x_{\tau}(t)) + f_1(x_{\tau}(t-\tau)) + u_{\tau}(t) f_2(x_{\tau}(t)) \medskip \\
x_{\tau}(t) = \phi(t) \ , \ t \in [-M,0] \ , \quad u_{\tau}(\cdot) \in L^2([0,T],\mathbb{R})
\end{cases}
\end{equation}
where $f_0$, $f_1$ and $f_2$ are smooth vector fields. This dynamics describes a great number of phenomena in physics, biology and economics \cite{erneux}, and it is widely used in engineering for modeling.

Fix a vector $\bar x \in \mathbb{R}^n$. We consider the Delayed Optimal Control Problem (\textbf{OCP})$_{\tau}$ of steering the control system (\ref{dynDelay}) to $\bar x$, minimizing the cost function
\begin{equation*} \label{costDelay}
C_T(u_{\tau}(\cdot)) = \int_0^T u^2_{\tau}(t) \; dt
\end{equation*}

This kind of control problems comes into play, for example, when the energy of a physical system is minimized or the value of the utility of an exhaustible resource is maximized.

The literature is abundant of numerical methods to solve (\textbf{OCP})$_{\tau}$, for example \cite{emmRev}, \cite{ex0}, \cite{ex1}, \cite{ex2}, \cite{ex3}, \cite{ex4}, \cite{ex5}. Nevertheless, some applications, like atmospheric reentry and satellite launching \cite{emmRev}, require great accuracy like \textit{indirect methods} can provide. It is then interesting understanding how to solve (\textbf{OCP})$_{\tau}$ via these procedures. Our impression is that the resolution of (\textbf{OCP})$_{\tau}$ by indirect methods has not yet been adequately addressed in the literature.

The basic version of indirect methods, the \textit{shooting method}, consists of solving, handling Newton-like algorithms, the two-point or multi-point boundary value problem arising from applying the \textit{Pontryagin Maximum Principle} (PMP) \cite{pontr}. Implementing such a routine on (\textbf{OCP})$_{\tau}$ originates two main troublesome drawbacks that make the classical indirect approach unusable.

The first one is the most critical and is induced by the PMP. It is known \cite{pontr} that, provided an optimum $(x_{\tau}(\cdot),u_{\tau}(\cdot))$, the delayed dual formulation ensures the existence of an absolutely continuous function $p_{\tau}(\cdot)$ of $[0,T]$ and a non-positive real number $p^0_{\tau}$ satisfying
\begin{eqnarray} \label{dynDual}
\begin{cases}
\displaystyle \dot{p}_{\tau}(t) = \displaystyle -\frac{\partial H}{\partial x}(x_{\tau}(t),x_{\tau}(t-\tau),u_{\tau}(t),p_{\tau}(t),p^0_{\tau}) \medskip \\
\displaystyle \hspace{15pt} -\frac{\partial H}{\partial y}(x_{\tau}(t+\tau),x_{\tau}(t),u_{\tau}(t+\tau),p_{\tau}(t+\tau),p^0_{\tau}) \ , \medskip \\
\hspace{15pt} t \in [0,T-\tau] \medskip \\
\displaystyle \dot{p}_{\tau}(t) = \displaystyle -\frac{\partial H}{\partial x}(x_{\tau}(t),x_{\tau}(t-\tau),u_{\tau}(t),p_{\tau}(t),p^0_{\tau}) \ , \medskip \\
\hspace{15pt} t \in (T-\tau, T]
\end{cases}
\end{eqnarray}
where $H(x,y,u,p,p^0) = \langle p , f_0(x) + f_1(y) + u f_2(x) \rangle + p^0 u$ is the Hamiltonian. Assuming that $u_{\tau}(\cdot)$ is known as a function of $x_{\tau}(\cdot)$ and $p_{\tau}(\cdot)$, each iteration of a shooting method consists of solving the coupled dynamics (\ref{dynDelay})-(\ref{dynDual}), where a value of $p_{\tau}(T)$ is provided. This means that one has to solve a system of \textit{mixed type differential equations}. The difficulty which comes up from dealing with these equations is the lack of global information which do not allow a local integration by usual iterative methods \cite{stoer}. Some techniques to solve mixed type differential equations were developed \cite{theo1}, \cite{theo2}, \cite{lunel}. However, these approaches are limited to linear systems.

The second flaw of using indirect methods, common to any kind of shooting algorithm, is that they are possibly hard to initialize because of their dependence on Newton methods \cite{stoer}. Indeed, even if we assume that we are able to solve (\ref{dynDelay})-(\ref{dynDual}), a good numerical guess of function $p_{\tau}(\cdot)$ must be provided first, to make the whole procedure converge.

Let the Non-Delayed Optimal Control Problem \textnormal{(\textbf{OCP})}$_0$ be
\begin{eqnarray*}
\begin{cases}
\displaystyle \min \ \ C_T(u(\cdot)) = \hspace{-0.1cm} \int_0^T u^2(t) \; dt \ , \quad u(\cdot) \in L^2([0,T],\mathbb{R}) \medskip \\
\dot{x}(t) = f_0(x(t)) + f_1(x(t)) + u(t) f_2(x(t)) \medskip \\
x(t) = \phi(t) \ , \quad t \in [-M,0] \ \ , \qquad x(T) = \bar x
\end{cases}
\end{eqnarray*}

In many situations, exploiting the non-delayed PMP mixed to other techniques \cite{emmRev}, one is able to initialize efficiently a shooting method on (\textbf{OCP})$_0$. Thus, it is legitimate to wonder if one may solve (\textbf{OCP})$_{\tau}$ by indirect methods starting an iterative procedure on (\textbf{OCP})$_0$ like \textit{homotopy methods} \cite{allgower}. This approach is a way to address the previous flaws: on one hand, we could solve (\ref{dynDelay})-(\ref{dynDual}) by using the previous global state solution and, on the other hand, the related global adjoint solution could be used to guess the value of $p_{\tau}(\cdot)$ of the next iteration.

Homotopy methods consist in deforming the problem into a simpler one that we are able to solve ((\textbf{OCP})$_0$ in our case) and then solving a sequence of shooting problems, step by step, to come back to the original problem. However, without appropriate assumptions, they can fail whenever, during the iteration path, bifurcation points, singularities or different connected components are encountered \cite{emmRev}, \cite{allgower}.

The idea proposed in this paper consists of introducing a new method coupled with a convergence result that allows to solve successfully (\textbf{OCP})$_{\tau}$ using indirect methods, by applying homotopy procedures which start from its non-delayed version (\textbf{OCP})$_0$.

The paper is organised as follows. In Section II the assumptions needed for the convergence theorem are provided followed by the main statement; then, the routine used to solve (\textbf{OCP})$_{\tau}$ is derived. Section III contains the proof of the main result. In Section IV the numerical effectiveness of the aforementioned procedure is showed by proposing a numerical example. Finally, Section V proposes some conclusions and perspectives while Section VI contains two technical results used in the main proof. 

\section{Assumptions, Main Result and Algorithm}

The main result that we present requires assumptions concerning (\textbf{OCP})$_0$ only. They are the following ones
\begin{itemize}
\item[$(H_1)$] $f_0$, $f_1$ and $f_2$ are bounded in $\mathbb{R}^n$;
\item[$(H_2)$] (\textbf{OCP})$_0$ has a unique solution, denoted $(x(\cdot),u(\cdot))$;
\item[$(H_3)$] The optimal trajectory $x(\cdot)$ has a unique extremal lift (up to a multiplicative scalar) defined in $[0,T]$, which is moreover normal, denoted $(x(\cdot), p(\cdot), -1, u(\cdot))$, solution of the Pontryagin Maximum Principle.
\end{itemize}

We want to stress the fact that we do not need any assumption concerning (\textbf{OCP})$_{\tau}$. Assumption $(H_1)$ plays only a technical role within the proof and it is not limiting: provided that all considered trajectories remain in some bounded set, one may multiply the dynamics by a vanishing smooth function. $(H_2)$, $(H_3)$ are usually stronger to satisfy and we suggest in Section III C how they can be weakened.

\begin{theorem} \label{mainTh}
Under Assumptions $(H_1)$, $(H_2)$ and $(H_3)$, there exists $\bar \tau > 0$ such that, for every $\tau \in (0,\bar \tau)$, (\textbf{OCP})$_{\tau}$ has at least one solution $(x_{\tau}(\cdot),u_{\tau}(\cdot))$, every extremal lift of which is normal. Let $(x_{\tau}(\cdot),u_{\tau}(\cdot),p_{\tau}(\cdot),-1)$ such a normal extremal lift. Then, as $\tau$ tends to $0$,
\begin{itemize}
\item $x_{\tau}(\cdot)$ converges uniformly to $x(\cdot)$ on $[-M,T]$;
\item $p_{\tau}(\cdot)$ converges uniformly to $p(\cdot)$ on $[0,T]$;
\item $u_{\tau}(\cdot)$ converges to $u(\cdot)$ in $L^2$ for the weak topology.
\end{itemize}
\end{theorem}

It is crucial to note that Theorem \ref{mainTh} gives the non-trivial conclusion that the adjoint vectors $p_{\tau}(\cdot)$ of (\textbf{OCP})$_{\tau}$ converge to the adjoint $p(\cdot)$ of the non-delayed version (\textbf{OCP})$_0$, independently from the sequence of $\tau$ chosen. Since, as we said before, there exist several procedures to initialize efficiently shooting methods on (\textbf{OCP})$_0$, the following algorithm arises straightforwardly

%\vspace{50pt}

\begin{enumerate}
\vspace{-7pt}
\begin{algorithm}
\item[1.] Set $k=0$, $\tau_k = 0$. Solve (\textbf{OCP})$_{\tau_0}$ by indirect methods and denote $p_{\tau_0}(\cdot)$ its adjoint vector solution;
\end{algorithm}
\vspace{-23pt}
\begin{algorithm}
\item[2.] Given a desired delay $\tau_d < \bar \tau$ to achieve, start a homotopy method on $\tau$ :

\While{$\tau_{k}<\tau_d$}{

\smallskip $\bullet$ Set $k = k+1$ and update $\tau_{k}$ ;

\medskip $\bullet$ Compute $p_{\tau_{k}}(\cdot)$ by solving (\textbf{OCP})$_{\tau_{k}}$ using a shooting method initialized with $p_{\tau_{k-1}}(\cdot)$.

}
\end{algorithm}
\end{enumerate}
\vspace{-12pt}
More precisely, within the shooting method at step 2. , the adjoint vector of the previous iteration is used as a \textit{global guess}, thanks to which, the integration of (\ref{dynDelay})-(\ref{dynDual}) is computed by replacing the time-forward elements by known terms.

The convergence of this method is ensured by Theorem~\ref{mainTh}, independently from how $\tau$ converges to zero. This allows to use every numerical homotopy approach to solve (\textbf{OCP})$_{\tau}$.

\section{Proof of Theorem \ref{mainTh}}

In this section, we focus on the proof of Theorem \ref{mainTh}. Two technical results used in the proof are recalled in the Appendix (Section VI).

\subsection{Existence Theorem}

The first step consists in establishing the existence of solutions of (\textbf{OCP})$_{\tau}$. The approach exploits a sensitivity analysis of the \textit{end-point mapping} as done in the proof of the PMP.

\begin{theorem} \label{existence}
Under Assumptions $(H_1)$, $(H_2)$ and $(H_3)$, there exist $\bar \tau > 0$ and $R > 0$ such that, for every $\tau \in (0,\bar \tau)$, (\textbf{OCP})$_{\tau}$ has at least one solution $(x_{\tau}(\cdot),u_{\tau}(\cdot))$ which is continuous with respect to $\tau$ and such that $\| u_{\tau}(\cdot) \|_{L^2} \leq R$.
\end{theorem}
\textbf{Proof:} Let $v(\cdot) \in L^2([0,T],\mathbb{R})$. For every $\tau \in [0,M]$, we denote $x_{\tau,v}(\cdot)$ the solution of
$$
\begin{cases}
\dot{y}(t) = f_0(y(t)) + f_1(y(t-\tau)) + v(t) f_2(y(t)) \medskip \\
y(t) = \phi(t) \ , \quad t \in [-M,0]
\end{cases}
$$

It is straightforward that, thanks to $(H_1)$, $x_{\tau,v}(\cdot)$ is defined on the whole interval $[-M,T]$. Then, we can define
$$
\tilde E_T : [0,M] \times L^2([0,T],\mathbb{R}) \rightarrow \mathbb{R}^n : (\tau,v(\cdot)) \mapsto x_{\tau,v}(T) - \bar x
$$
Thanks to Theorem \ref{sd} (Section VI), $\tilde E_T$ is continuous and continuously differentiable with respect to $v$, for every $\tau~\in~[0,M]$.

By $(H_3)$, $(x(\cdot),u(\cdot),p(\cdot),-1)$ is the unique extremal lift of (\textbf{OCP})$_0$, which is moreover normal. From the PMP, it follows that, for every $v(\cdot) \in L^2([0,T],\mathbb{R})$
$$
p(T).\frac{\partial \tilde E_T}{\partial v}(0,u(\cdot)) \cdot v(\cdot) - dC_T(u(\cdot)) \cdot v(\cdot) = 0
$$

If one supposes that $\frac{\partial \tilde E_T}{\partial v}(0,u(\cdot))$ is not surjective, it would exist a row vector $\psi \in \mathbb{R}^n \setminus \{ 0 \}$ such that
$$
\psi.\frac{\partial \tilde E_T}{\partial v}(0,u(\cdot)) \cdot v(\cdot) = 0
$$

But, this implies that $(x(\cdot),u(\cdot))$ would have an abnormal extremal lift, which raises a contradiction. It follows
$$
\displaystyle \frac{\partial \tilde E_T}{\partial v}(0,u(\cdot)) \cdot L^2([0,T],\mathbb{R}) = \mathbb{R}^n
$$

Hence, from Theorem \ref{ift} (Section VI), it exists a real number $\bar \tau > 0$ and a $C^0$ function $\tau \in [0,\bar \tau) \mapsto u_{\tau}(\cdot) \in L^2([0,T],\mathbb{R})$ such that, for every $\tau \in [0,\bar \tau)$, it holds $x_{\tau,u_{\tau}}(T) = \bar x$. Moreover, by possibly taking a smaller $\bar \tau$, we assume that
$$
\| u_{\tau}(\cdot) \|_{L^2} \leq \| u(\cdot) \|_{L^2} + 1 = R \ , \ \forall \tau \in (0,\bar \tau)
$$

Consider now the following problem
$$
(P^R_{\tau})
\begin{cases}
\displaystyle \min \int_0^T v^2(t) \; dt \ , \ \; v(\cdot) \in L^2([0,T],\mathbb{R}) \ , \ \; \| v(\cdot) \|_{L^2} \leq R \medskip \\
\dot{y}(t) = f_0(y(t)) + f_1(y(t-\tau)) + v(t) f_2(y(t)) \medskip \\
y(t) = \phi(t) \ , \quad t \in [-M,0] \ \ , \qquad y(T) = \bar x
\end{cases}
$$

The previous argument shows that, for every $\tau~\in~(0,\bar \tau)$, $(P^R_{\tau})$ is controllable. Moreover, thanks to $(H_1)$, all admissible trajectories of $(P^R_{\tau})$ are uniformly bounded.

Take a sequence of controls $(v_k(\cdot))_{k \in \mathbb{N}}~\subseteq~L^2([0,T],\mathbb{R})$ admissible for $(P^R_{\tau})$ such that their costs converge to the infimum taken on the set of admissible controls of $(P^R_{\tau})$. Let $ (y_k(\cdot))_{k \in \mathbb{N}}$ the sequence of the associated solutions. Since $(v_k(\cdot))_{k \in \mathbb{N}}$ is bounded in $L^2$, it exists a function $\bar v(\cdot)~\in~L^2([0,T],\mathbb{R})$ such that, up to a sequence, $v_k(\cdot)$ converges to $\bar v(\cdot)$ in $L^2$ for the weak topology. Clearly
\begin{equation} \label{limit0}
y_k(t) = \phi(t) \cdot \mathds{1}_{[-M,0]}(t) + \mathds{1}_{(0,T]}(t) \cdot \bigg[ \phi(0)
\end{equation}
$$
+ \int^t_0 \bigg(f_0(y_k(s)) + f_1(y_k(s-\tau)) + v_k(s) f_2(y_k(s)) \bigg) \ ds \bigg]
$$

Since the sequence $(\dot{y}_k(\cdot))_{k \in \mathbb{N}}$ is bounded in $L^2([0,T],\mathbb{R})$, $(y_k(\cdot))_{k \in \mathbb{N}}$ is bounded in $H^1([0,T],\mathbb{R}^n)$. Hence, it exists a function $\bar y(\cdot) \in H^1([0,T],\mathbb{R}^n)$ such that, up to a sequence, $y_k(\cdot)$ converges to $\bar y(\cdot)$ in $H^1$ for the weak topology. The embedding of $H^1([0,T],\mathbb{R}^n)$ into $C^0([0,T],\mathbb{R}^n)$ is compact. Then, it follows that $y_k(\cdot)$ converges uniformly to $\bar y(\cdot)$ within $[0,T]$. In particular $\bar y(0) = \phi(0)$, $\bar y(T) = \bar x$. Hence, we can extend continuously $\bar y(\cdot)$ to $[-M,T]$ setting $\bar y(t) = \phi(t)$ for $t \in [-M,0]$. Thanks to the definition of weak convergence, passing to the limit within (\ref{limit0}) gives easily
$$
\bar y(t) = \phi(t) \cdot \mathds{1}_{[-M,0]}(t) + \mathds{1}_{(0,T]}(t) \cdot \bigg[ \phi(0)
$$
$$
+ \int^t_0 \bigg(f_0(\bar y(s)) + f_1(\bar y(s-\tau)) + \bar v(s) f_2(\bar y(s)) \bigg) \ ds \bigg]
$$

Then, $(\bar y(\cdot), \bar v(\cdot))$ is an optimal solution of $(P^R_{\tau})$. Since an optimal solution of $(P^R_{\tau})$ is also optimal for (\textbf{OCP})$_{\tau}$, the conclusion follows.~$_{\Box}$

\subsection{Convergence Theorem}

Once the existence of solutions of (\textbf{OCP})$_{\tau}$ is established, we can use their $C^0$ dependence with respect to $\tau$ to prove their convergence to solutions of (\textbf{OCP})$_0$. \\

\noindent \textbf{Proof of Theorem \ref{mainTh}:} We divide the proof in three parts. \\

\noindent \textbf{Existence:} Thanks to Theorem \ref{existence}, there exist $\bar \tau > 0$, $R > 0$ such that, for every $\tau \in (0,\bar \tau)$, (\textbf{OCP})$_{\tau}$ has at least one solution $(x_{\tau}(\cdot),u_{\tau}(\cdot))$, for which $\| u_{\tau}(\cdot) \|_{L^2} \leq R$. \\

\noindent \textbf{Convergence of $\bm{(x_{\tau}(\cdot),u_{\tau}(\cdot))}$:} Let $(\tau_k)_{k \in \mathbb{N}}$ be a sequence of positive real numbers converging to 0. We denote $(x_{\tau_k}(\cdot),u_{\tau_k}(\cdot))_{k \in \mathbb{N}}$ the associated sequence of optimal solutions of (\textbf{OCP})$_{\tau_k}$. Since $(u_{\tau_k}(\cdot))_{k \in \mathbb{N}}$ is bounded in $L^2$, it exists a function $\bar u(\cdot) \in L^2([0,T],\mathbb{R})$ such that, up to a subsequence, $u_{\tau_k}(\cdot)$ converges to $\bar u(\cdot)$ in $L^2$ for the weak topology. Thanks to $(H_1)$ and $\| u_{\tau}(\cdot) \|_{L^2} \leq R$, from
\begin{equation*} \label{limit}
x_{\tau_k}(t) = \phi(t) \cdot \mathds{1}_{[-M,0]}(t) + \mathds{1}_{(0,T]}(t) \cdot \bigg[ \phi(0)
\end{equation*}
$$
+ \int^t_0 \bigg(f_0(x_{\tau_k}(s)) + f_1(x_{\tau_k}(s-\tau_k)) + u_{\tau_k}(s) f_2(x_{\tau_k}(s)) \bigg) \ ds \bigg]
$$
we obtain that $(x_{\tau_k}(\cdot))_{k \in \mathbb{N}}$ is bounded in $C^0([-M,T],\mathbb{R}^n)$ and consequently in $L^2([0,T],\mathbb{R}^n)$. Moreover, $(\dot{x}_{\tau_k}(\cdot))_{k \in \mathbb{N}}$ has an upper bound in $L^2([0,T],\mathbb{R}^n)$. Then, $(x_{\tau_k}(\cdot))_{k \in \mathbb{N}}$ is bounded in $H^1([0,T],\mathbb{R}^n)$. This implies that it exists a function $\bar x~\in~H^1([0,T],\mathbb{R}^n)$ such that, up to a subsequence, $x_{\tau_k}(\cdot)$ converges to $\bar x(\cdot)$ in $H^1$ for the weak topology. Since the embedding of $H^1([0,T],\mathbb{R}^n)$ into $C^0([0,T],\mathbb{R}^n)$ is compact, it follows that $x_{\tau_k}(\cdot)$ converges uniformly to $\bar x(\cdot)$ within $[0,T]$. In particular $\bar x(0) = \phi(0)$, $\bar x(T) = \bar x$. Hence, we can extend continuously $\bar x(\cdot)$ to $[-M,T]$ setting $\bar x(t) = \phi(t)$ for $t \in [-M,0]$. From this, we have, for every $t \in [0,T]$
$$
\bigg| \displaystyle \lim_{\tau_k \rightarrow 0} \int^t_0 f_1(x_{\tau_k}(s-\tau_k)) \ ds - \int^t_0 f_1(\bar x(s)) \ ds \bigg|
$$
$$
\leq \displaystyle \lim_{\tau_k \rightarrow 0} \int^t_0 | f_1(x_{\tau_k}(s-\tau_k)) - f_1(\bar x(s-\tau_k)) | \ ds
$$
$$
+ \lim_{\tau_k \rightarrow 0} \int^t_0 | f_1(\bar x(s-\tau_k)) - f_1(\bar x(s)) | \ ds
$$
$$
\leq \displaystyle \lim_{\tau_k \rightarrow 0} T \| f_1(x_{\tau_k}(\cdot)) - f_1(\bar x(\cdot)) \|_{C^0} + 0 = 0
$$
by the dominated convergence theorem, and
$$
\bigg| \displaystyle \lim_{\tau_k \rightarrow 0} \int^t_0 u_{\tau_k}(s) f_2(x_{\tau_k}(s)) \ ds - \int^t_0 \bar u(s) f_2(\bar x(s)) \ ds \bigg|
$$
$$
\leq \displaystyle \lim_{\tau_k \rightarrow 0} \int^t_0 | u_{\tau_k}(s) f_2(x_{\tau_k}(s)) - u_{\tau_k}(s) f_2(\bar x(s)) | \ ds
$$
$$
+ \bigg| \lim_{\tau_k \rightarrow 0} \int^t_0 u_{\tau_k}(s) f_2(\bar x(s)) \ ds - \int^t_0 \bar u(s) f_2(\bar x(s)) \ ds \bigg|
$$
$$
\leq \displaystyle \lim_{\tau_k \rightarrow 0} \sqrt{T} R \| f_2(x_{\tau_k}(\cdot)) - f_2(\bar x(\cdot)) \|_{C^0} + 0 = 0
$$
by the definition of weak convergence. Then
$$
\bar x(t) = \phi(t) \cdot \mathds{1}_{[-M,0]}(t) + \mathds{1}_{(0,T]}(t) \cdot \bigg[ \phi(0)
$$
$$
+ \int^t_0 \bigg(f_0(\bar x(s)) + f_1(\bar x(s)) + \bar u(s) f_2(\bar x(s)) \bigg) \ ds \bigg]
$$
when $t \in [-M,T]$. Hence, $\bar x(\cdot)$ is an admissible trajectory of (\textbf{OCP})$_0$, generated by control $\bar u(\cdot)$. Now, from Theorem \ref{existence}, for every $k \in \mathbb{N}$, it exists $v_{\tau_k}(\cdot) \in L^2([0,T],\mathbb{R})$ admissible for (\textbf{OCP})$_{\tau_k}$ such that $\| v_{\tau_k}(\cdot) \|_{L^2}$ converges to $\| u(\cdot) \|_{L^2}$. Since
$$
\| \bar u(\cdot) \|_{L^2} \leq \liminf \| u_{\tau_k}(\cdot) \|_{L^2} \leq \| v(\cdot) \|_{L^2} \ ,
$$
$$
\forall \ k \in \mathbb{N} \ , \ v(\cdot) \ \ \textnormal{admissible for} \ \textnormal{(\textbf{OCP})$_{\tau_k}$}
$$
passing to the limit, we infer that $\| \bar u(\cdot) \|_{L^2} \leq \| u(\cdot) \|_{L^2}$. Hence, $(H_2)$ gives $\bar x(\cdot) = x(\cdot)$ in $[-M,T]$ and $\bar u(\cdot) = u(\cdot)$ almost everywhere in $[0,T]$. Since $(\tau_k)_{k \in \mathbb{N}}$ is an arbitrary sequence, these convergence results hold for the whole family $(x_{\tau}(\cdot),u_{\tau}(\cdot))$. \\

\noindent \textbf{Convergence of $\bm{p_{\tau}(\cdot)}$:} From formulation (\ref{dynDual}) applied to (\textbf{OCP})$_{\tau}$, we know that, for every $\tau \in (0,\bar \tau)$, there exists an extremal lift $(x_{\tau}(\cdot),u_{\tau}(\cdot),p_{\tau}(\cdot),p^0_{\tau})$ satisfying
\begin{eqnarray*}
(D_{\tau})
\begin{cases}
\displaystyle \dot{p}_{\tau}(t) = \displaystyle -\bigg\langle p_{\tau}(t) , \frac{\partial f_0}{\partial x}(x_{\tau}(t)) + u_{\tau}(t) \frac{\partial f_2}{\partial x}(x_{\tau}(t)) \bigg\rangle \medskip \\
\displaystyle \hspace{1.2cm} -\bigg\langle p_{\tau}(t+\tau) , \frac{\partial f_1}{\partial x}(x_{\tau}(t)) \bigg\rangle \ , \quad t \in [0,T-\tau] \medskip \\
\displaystyle \dot{p}_{\tau}(t) = \displaystyle -\bigg\langle p_{\tau}(t) , \frac{\partial f_0}{\partial x}(x_{\tau}(t)) + u_{\tau}(t) \frac{\partial f_2}{\partial x}(x_{\tau}(t)) \bigg\rangle \ , \medskip \\
\displaystyle \hspace{1.45cm} t \in (T-\tau, T]
\end{cases}
\end{eqnarray*}
%\begin{eqnarray*}
%(P^A_{\tau})
%\begin{cases}
%\displaystyle \dot{p}_{\tau}(t) = \displaystyle -\frac{\partial H}{\partial x}(x_{\tau}(t),x_{\tau}(t-\tau),u_{\tau}(t),p_{\tau}(t),p^0_{\tau}) \medskip \\
%\displaystyle \hspace{1.45cm} -\frac{\partial H}{\partial y}(x_{\tau}(t+\tau),x_{\tau}(t),u_{\tau}(t+\tau),p_{\tau}(t+\tau),p^0_{\tau}) \ , \medskip \\
%\displaystyle \hspace{1.45cm} t \in [0,T-\tau] \medskip \\
%\displaystyle \dot{p}_{\tau}(t) = \displaystyle -\frac{\partial H}{\partial x}(x_{\tau}(t),x_{\tau}(t-\tau),u_{\tau}(t),p_{\tau}(t),p^0_{\tau}) \ , \medskip \\
%\displaystyle \hspace{1.45cm} t \in [T-\tau, T]
%\end{cases}
%\end{eqnarray*}
%where $H(x(t),y(t),u(t),p(t),p^0) = \langle p(t) , f_0(x(t)) + f_1(y(t)) + u(t) f_2(x(t)) \rangle + p^0 C_T(u(t))$.

In the following, we denote
$$
E_T(\tau,v(\cdot)) : [0,\bar \tau) \times L^2([0,T],\mathbb{R}) \rightarrow \mathbb{R}^n : (\tau,v(\cdot)) \mapsto x_{\tau,v}(T)
$$
where $x_{\tau,v}(\cdot)$ is the trajectory associated to $v(\cdot)$ satisfying (\ref{dynDelay}). By Theorem \ref{sd}, $E_T(\tau,v(\cdot))$ is continuous in both variables and has continuous derivative with respect to $v$. Moreover, for every sequence $(u_{k}(\cdot))_{k \in \mathbb{N}}$ of controls weakly convergent to a control $u(\cdot)$ and every sequence of real numbers $(\tau_k)_{k \in \mathbb{N}}~\subseteq~[0,M]$ converging to $\tau~\in~[0,M]$
$$
\frac{\partial E_T}{\partial v}(\tau_k,u_k(\cdot)) \longrightarrow \frac{\partial E_T}{\partial v}(\tau,u(\cdot))
$$

For sake of clarity, we split this proof in subparts. \\

1) \textit{Every extremal lift $(x_{\tau}(\cdot),u_{\tau}(\cdot),p_{\tau}(\cdot),p^0_{\tau})$ of any solution of (\textbf{OCP})$_{\tau}$ is normal, i.e. $p^0_{\tau} = -1$.}

We argue by contradiction. Suppose that for every integer $k$ there exists $\tau_k \in (0,1/k)$ such that $(x_{\tau_k}(\cdot),u_{\tau_k}(\cdot),p_{\tau_k}(\cdot),0)$ is an extremal lift associated to the solution $(x_{\tau_k}(\cdot),u_{\tau_k}(\cdot))$ of (\textbf{OCP})$_{\tau_k}$. Then, from the PMP, we have
\begin{equation} \label{exit}
\forall v(\cdot) \in L^2([0,T],\mathbb{R}) \ , \quad p_{\tau_k}(T).\frac{\partial E_T}{\partial v}(\tau_k,u_{\tau_k}(\cdot)) \cdot v(\cdot) = 0
\end{equation}
We denote $\psi_{\tau_k} = p_{\tau_k}(T)$. Since $p_{\tau_k}(T)$ is unique up to a multiplicative constant, we can suppose $\| \psi_{\tau_k} \| = 1$ for every integer $k$. Then, up to a subsequence, $(\psi_{\tau_k})_{\mathbb{N}}$ converges to some vector $\psi \neq 0$. Passing to the limit within (\ref{exit}), thanks to the previous results, it follows immediately
$$
\forall v(\cdot) \in L^2([0,T],\mathbb{R}) \ , \quad \psi.\frac{\partial E_T}{\partial v}(0,u(\cdot)) \cdot v(\cdot) = 0
$$
which raises a contradiction because of $(H_3)$. \\

We remark that, thanks to the argument above, from the delayed version of the PMP we have that, for every $\tau~\in~(0,\bar \tau)$ and every $v(\cdot) \in L^2([0,T],\mathbb{R})$, it holds
\begin{equation} \label{exitComp}
p_{\tau}(T).\frac{\partial E_T}{\partial v}(\tau,u_{\tau}(\cdot)) \cdot v(\cdot) - dC_T(u_{\tau}(\cdot)) \cdot v(\cdot) = 0
\end{equation}

\vspace{5pt}

2) \textit{Setting $p^0_{\tau} = -1$, the set of all possible $p_{\tau}(T)$, with $\tau \in (0,\bar \tau)$, is bounded.}

We proceed again by contradiction. Assume that there exists a sequence $(\tau_k)_{\mathbb{N}}$ of positive real numbers converging to 0 such that $\| p_{\tau_k}(T) \|$ tends to $+\infty$. Since the sequence $\bigg( \frac{p_{\tau_k}(T)}{\| p_{\tau_k}(T) \|} \bigg)_{k \in \mathbb{N}}$ is bounded in $\mathbb{R}^n$, up to a subsequence, it converges to some unit vector $\psi$. Then, dividing by $\| p_{\tau_k}(T) \|$ and passing to the limit within (\ref{exitComp}), thanks to the previous results, we have
$$
\forall v(\cdot) \in L^2([0,T],\mathbb{R}) \ , \quad \psi.\frac{\partial E_T}{\partial v}(0,u(\cdot)) \cdot v(\cdot) = 0
$$
which raises a contradiction because of $(H_3)$. \\

3) \textit{$p_{\tau}(\cdot)$ converges uniformly to $p(\cdot)$ on $[0,T]$, where $(x(\cdot),u(\cdot),p(\cdot),-1)$ is the unique extremal lift of $(x(\cdot),u(\cdot))$.}

For every $\tau \in (0,\bar \tau)$, set $\psi_{\tau} = p_{\tau}(T)$. We know that the adjoint system of (\textbf{OCP})$_{\tau}$ is represented by $(D_{\tau})$.

From the previous argumentation, the family of all $\psi_{\tau}$, $0~<~\tau~<~\bar \tau$ is bounded. Let $\psi$ be a closure point of that family, and $(\tau_k)_{k \in \mathbb{N}}$ a sequence of positive real numbers converging to 0 such that $\psi_{\tau_k}$ tends to $\psi$. Using the continuous dependence from initial data, we infer that the sequence $(p_{\tau_k}(\cdot))_{k \in \mathbb{N}}$ converges uniformly to the solution $z(\cdot)$ of the Cauchy problem
$$
\begin{cases}
\displaystyle \dot{z}(t) = \displaystyle -\bigg\langle z(t) , \frac{\partial f_0}{\partial x}(x(t)) + \frac{\partial f_1}{\partial x}(x(t)) + u(t) \frac{\partial f_2}{\partial x}(x(t)) \bigg\rangle \medskip \\
z(T) = \psi
\end{cases}
$$
Moreover, passing to the limit within (\ref{exitComp}) gives
$$
\psi.\frac{\partial E_T}{\partial v}(0,u(\cdot)) \cdot v(\cdot) - dC_T(u(\cdot)) \cdot v(\cdot) = 0
$$
for every $v(\cdot) \in L^2([0,T],\mathbb{R})$. It follows immediately that $(x(\cdot),u(\cdot),z(\cdot),-1)$ is a normal extremal lift of (\textbf{OCP})$_0$. Exploiting $(H_3)$, we obtain $z(\cdot) = p(\cdot)$ in $[0,T]$. \\

The conclusion follows. $_\Box$

\subsection{Remarks}

Assumptions $(H_2)$, $(H_3)$ allow to identify the limit of the sequence of solutions of (\textbf{OCP})$_{\tau}$. However, they are usually not satisfied or difficult to prove. Nevertheless, one can prove that these assumptions remain satisfied in certain general frameworks. Moreover, they can be weakened to modify Theorem \ref{mainTh} in terms of closure points. We refer to remarks 2.8, 2.19 of \cite{Emm} and to remarks 1.4, 2.6 of \cite{EmmSilva} for details.

\section{Numerical Example}

In this section we apply the method introduced by Theorem \ref{mainTh} to a two dimensional rendezvous problem, usually encountered in aerospace and robotic applications.

The framework is as follows. Consider a vehicle which moves with constant velocity $v_0$ and whose path lies in a plane $(x,y)$. The mission consists in reaching a point $(\bar x,\bar y)$ with a specified orientation $\bar \theta$ starting from $(x_0,y_0,\theta_0)$ in a fixed time $T>0$. When the actuators change set-up, mechanical frictions or aerodynamic effects prevent the instantaneous displacement of the vehicle; this is modeled by a constant delay $\tau$ in changing orientations. The dynamics is
\begin{eqnarray} \label{dynEx}
\begin{cases}
\dot{x}_{\tau}(t) = v_0 \cos(\theta_{\tau}(t)) \ , \quad \dot{y}_{\tau}(t) = v_0 \sin(\theta_{\tau}(t)) \medskip \\
\dot{\theta}_{\tau}(t) = c_0 v_0 \delta_{\tau} (t-\tau) \ , \quad \dot{\delta}_{\tau}(t) = u_{\tau}(t)
\end{cases}
\end{eqnarray}
where $c_0>0$ is constant (that we set equal to one, without loss of generality) and $\delta_{\tau}$ is the steering angle.

One can ask to minimize the energy used to change path; this is equivalent to solve the problem of minimizing the integral of the square of $u_{\tau}(\cdot)$ subject to dynamics (\ref{dynEx}). Then, the algorithm introduced in Section II can be applied to solve this specific (\textbf{OCP})$_{\tau}$. In particular, without loss of generality (see Section II), a linear continuation is used.

We impose as initial and final states respectively $(x_0,y_0,\theta_0,\delta_0) = (0,0,\pi/4,5 \cdot 10^{-4})$, constant in $[-\tau,0]$, and $(\bar x,\bar y,\bar \theta,\bar \delta) = (1500,1000,\pi/20,0)$, with $T = 19 \ s$. We run three tests for which we suppose a constant velocity of $100 \ m/s$ and delays respectively of $\tau~=~0 \ s$, $\tau = 2 \ s$, $\tau = 4 \ s$.

The optimal values of the costs obtained are $C_T(u_0(\cdot)) = 6.09179 \cdot 10^{-7}$, $C_T(u_{2}(\cdot)) = 6.05183 \cdot 10^{-7}$ and $C_T(u_{4}(\cdot)) = 8.68821 \cdot 10^{-7}$ while the optimal values of the initial adjoint vectors obtained are
\begin{center}
\begin{tabular}{@{}c|@{}c|@{}c|@{}c|@{}c}
$\tau$ ($s$) & $p_x(0) \cdot 10^{8}$ & $p_y(0) \cdot 10^{8}$ & $p_{\theta}(0) \cdot 10^{6}$ & $p_{\delta}(0) \cdot 10^{4}$ \\
\hline
0 & $-1.4576$ & $-1.0064$ & $2.8854$ & $2.9331$ \\
2 & $-1.4406$ & $-0.95041$ & $3.1203$ & $-1.2453$ \\
4 & $0.30418$ & $-0.24459$ & $-5.2204$ & $-13.948$
\end{tabular}
\end{center}

Figures 1 shows the optimal controls. It is clear that, adding delay, the solutions drift apart from the non-delayed one through a \textit{non-intuitive way}, difficult to foresee. The approach proposed happens to be very efficient. For example, testing the previous case, a fast rate of convergence is ensured when $\tau \in (0,7]$, even if such extended range of values is physically meaningless.

\section{Conclusion and Perspectives}

In this paper, we have presented a new approach that allows to solve optimal control problems with delayed control-affine dynamics exploiting indirect methods coupled with homotopy procedures. The idea is to introduce the delay step by step with a continuation parameter, starting from the solution of the non-delayed problem. This procedure is suggested by a convergence result which allows to initialize successfully the related shooting method by providing a global guess of the adjoint vector.

Future works on this subject focus on extending the aforementioned procedure to general non-linear delayed optimal control problems with control constraints. The main step consists of adapting Theorem \ref{mainTh} to the control constraints framework, exploiting the \textit{needle-like variation} technique.

\section{Appendix}

Here, some technical results used in Section III are presented. We recall first the Implicit Function Theorem.
\begin{theorem} \label{ift}
Let $U$ be a Banach space and consider a mapping
$$
F : [0,M] \times U \rightarrow \mathbb{R}^n : (\varepsilon,x) \mapsto F(\varepsilon,x)
$$
satisfying the following assumptions:
\begin{itemize}
\item $F(0,x_0) = 0$;
\item $F$ is continuous;
\item For every $\varepsilon \in [0,M]$, $F$ is continuously differentiable with respect to $x$ at $x_0$ and $\displaystyle \frac{\partial F}{\partial x}(\varepsilon,x_0)$ is continuous with respect to $\varepsilon$;
\item $\displaystyle \frac{\partial F}{\partial x}(0,x_0)(U) = \mathbb{R}^n$.
\end{itemize}
Then, there exists $\varepsilon_0 > 0$ and a continuous function $g~:~[0,\varepsilon_0) \rightarrow U$ such that $g(0) = x_0$ and
$$
F(\varepsilon,g(\varepsilon)) = 0 \ , \ \forall \varepsilon \in [0,\varepsilon_0)
$$
\end{theorem}

We have the following useful result
\begin{theorem} \label{sd}
Fix $M>0$, $x_1 \in \mathbb{R}^n$ and consider, for every control $u(\cdot) \in L^2([0,T],\mathbb{R})$, the following problem
$$
(P)
\begin{cases}
\dot{x}_{\tau,u}(t) = f_0(x_{\tau,u}(t)) + f_1(x_{\tau,u}(t-\tau)) + u(t) f_2(x_{\tau,u}(t)) \medskip \\
x_{\tau,u}(t) = \phi(t) \ , \quad t \in [-M,0]
\end{cases}
$$
where $\tau \in [0,M]$ and $\phi$ is continuous in $[-M,0]$. If $f_0$, $f_1$ and $f_2$ are smooth and bounded, function
$$
E_T : [0,M] \times L^2([0,T],\mathbb{R}) \rightarrow \mathbb{R}^n : (\tau,v) \mapsto x_{\tau,v}(T)
$$
is well defined and continuous by the continuous dependence from initial data of $(P)$. Moreover, for every $\tau \in [0,M]$ and $u(\cdot) \in L^2([0,T],\mathbb{R})$, $E_T$ is continuously differentiable with respect to $v$ at $u(\cdot)$ and it holds
$$
\displaystyle \frac{\partial E_T}{\partial v}(\tau,u) \cdot v(\cdot) = \int_0^T X_u(T,s) f_2(x_{\tau,u}(s)) v(s) \ ds
$$
where $X_u(t,s)$ is the matrix solution of
$$
(P_{\textnormal{Lin}})
\begin{cases}
\dot{X}_u(t,s) = A^1_u(t) X_u(t,s) + A^2_u(t) X_u(t-\tau,s) \medskip \\
X_u(t,s) = \begin{cases} I \ , \quad t=s \medskip \\ 0 \ , \quad t<s \end{cases}
\end{cases}
$$
with $A^1_u(t) = df_0(x_{\tau,u}(t)) + u(t) df_2(x_{\tau,u}(t))$, $A^2_u(t)~=~df_1(x_{\tau,u}(t-\tau))$. Finally, $\displaystyle \frac{\partial E_T}{\partial v}(\tau,u)$ is continuous because of the continuous dependence of $X_u(t,s)$ by initial data. In particular, for every sequence $(u_{k}(\cdot))_{k \in \mathbb{N}}$ of controls weakly convergent to a control $u(\cdot)$ and every sequence of real numbers $(\tau_k)_{k \in \mathbb{N}} \subseteq [0,M]$ converging to $\tau~\in~[0,M]$, it holds
$$
\frac{\partial E_T}{\partial v}(\tau_k,u_k(\cdot)) \longrightarrow \frac{\partial E_T}{\partial v}(\tau,u(\cdot))
$$
\end{theorem}
\textbf{Proof:} The first part can be easily achieved following the procedure developed to prove Proposition 3.6 of \cite{artEmm}. \\

In order to prove the last result, we may remark that, adapting easily the proof of Proposition 3.3 of \cite{artEmm}, it follows that if $(u_k)_{k \in \mathbb{N}}~\subseteq~L^2([0,T],\mathbb{R})$ is a sequence of controls weakly converging to some control $u \in L^2([0,T],\mathbb{R})$ and $(\tau_k)_{k \in \mathbb{N}}~\subseteq~[0,M]$ is a sequence converging to $\tau~\in~[0,M]$, $x_{\tau_k,u_k}(\cdot)$ converges uniformly to $x_{\tau,u}(\cdot)$. Hence, it suffices to show that $X_{\tau_k,u_k}(T,\cdot)$ converges uniformly to $X_{\tau,u}(T,\cdot)$ for every sequence of controls $(u_k)_{k \in \mathbb{N}}~\subseteq~L^2([0,T],\mathbb{R})$ weakly converging to some control $u \in L^2([0,T],\mathbb{R})$ and every sequence $(\tau_k)_{k \in \mathbb{N}}~\subseteq~[0,M]$ converging to $\tau~\in~[0,M]$. Since $A^1_{\tau_k,u_k}(\cdot)$, $A^2_{\tau_k,u_k}(\cdot)$ converge uniformly respectively to $A^1_{\tau,u}(\cdot)$, $A^2_{\tau,u}(\cdot)$ (hence, they are bounded) and, by the Gronwall inequality, $X_{\tau_k,u_k}(\cdot,s)$ are uniformly bounded on $[s,T]$, for every $\varepsilon > 0$ it exists an integer $N_{\varepsilon}$ such that, for every $k \geq N_{\varepsilon}$ and $s \in [0,T]$, we have
$$
\| X_{\tau,u}(T,s) - X_{\tau_k,u_k}(T,s) \|
$$
$$
\leq \int^T_s \| A^1_{\tau_k,u_k}(t) (X_{\tau,u}(t,s) - X_{\tau_k,u_k}(t,s)) \| \ dt
$$
$$
+ \int^T_s \| A^2_{\tau_k,u_k}(t) (X_{\tau,u}(t-\tau,s) - X_{\tau_k,u_k}(t-\tau_k,s)) \| \ dt
$$
$$
+ \int^T_s \| (A^1_{\tau,u}(t) - A^1_{\tau_k,u_k}(t)) X_{\tau,u}(t,s) \| \ dt
$$
$$
+ \int^T_s \| (A^2_{\tau,u}(t) - A^2_{\tau_k,u_k}(t)) X_{\tau,u}(t-\tau,s) \| \ dt
$$
$$
\leq C_1 \int^T_s \| X_{\tau,u}(t,s) - X_{\tau_k,u_k}(t,s) \| \ dt + \varepsilon
$$
$$
+ C_2 \bigg( \int^s_{s-\tau} \| X_{\tau,u}(r,s) - X_{\tau_k,u_k}(r,s) \| \ dr
$$
$$
+ \int^{T-\tau}_s \| X_{\tau,u}(r,s) - X_{\tau_k,u_k}(r,s) \| \ dr \bigg)
$$
$$
\leq C \int^T_s \| X_{\tau,u}(t,s) - X_{\tau_k,u_k}(t,s) \| \ dt + \varepsilon
$$
The thesis follows from the Gronwall inequality. $_\Box$

\vspace{-10pt}
%\begin{minipage}[c]{1.\linewidth}
%\begin{figure}[H]
%\includegraphics[width=0.85\textwidth]{z.eps}
%\caption{Trajectory $(x_{\tau}(t),y_{\tau}(t))$. The empty-dot line represents the non-delayed solution, while the dashed and the solid lines represent the solutions with delays respectively of $\tau = 2 \ s$ and $\tau = 4 \ s$. The same legend holds for Figure 2.}
%\end{figure}
%\end{minipage}

\begin{minipage}[c]{1.\linewidth}
\begin{figure}[H]
\centering
\includegraphics[width=1.\textwidth]{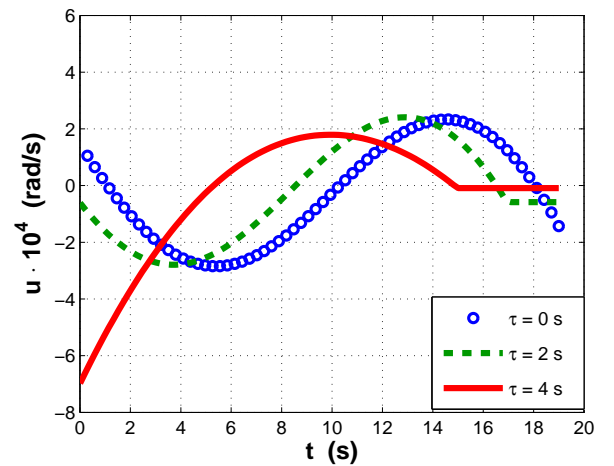}
\caption{Control $u_{\tau}(t)$.}
%\caption{\textcolor{red}{Control $u_{\tau}(t)$. The empty-dot line represents the non-delayed solution, while the dashed and the solid lines represent the solutions with delays respectively of $\tau = 2 \ s$ and $\tau = 4 \ s$.}}
\end{figure}
\end{minipage}


\begin{thebibliography}{99}
\bibitem{erneux} M. Malek-Zavarei, M. Jamshidi. \textit{Time Delay Systems: Analysis, Optimization and Applications}. (North-Holland Systems and Control Series), Elsevier Science, New York, 1987.

\bibitem{emmRev} E. Tr\'{e}lat. \textit{Optimal Control and Applications to Aerospace: some Results and Challenges}. Journal of Optimization Theory and Applications, September 2012, Volume 154, Issue 3, pp 713-758.

\bibitem{ex0} H.T. Banks. \textit{Approximation of nonlinear functional differential equation control systems}. Journal of Optimization Theory and Applications 29 (3) (1979) 383-408.

\bibitem{ex1} K. Kaji, K.H. Wong. \textit{Nonlinearly constrained time-delayed optimal control problems}. Journal of Optimization Theory and Applications 82 (2) (1994) 295–313.

\bibitem{ex2} C.L. Chen, D.Y. Su, C.Y. Chang. \textit{Numerical solution of time-delayed optimal control problems by iterative dynamic programming}. Optimal Control Applications and Methods 21 (2000) 91-105.

\bibitem{ex3} G.N. Elnagar, M.A. Kazemi. \textit{Numerical solution of time-delayed functional differential equation control systems}. Journal of Computational and Applied Mathematics 130 (2001) 75-90.

\bibitem{ex4} H.W.J. Lee, K.H. Wong. \textit{Semi-infinite programming approach to nonlinear time-delayed optimal control problems with linear continuous constraints}. Optimization Methods and Software 21 (5) (2006).

\bibitem{ex5} L. G{\"o}llmann, D. Kern, H. Maurer. \textit{Optimal control problems with delays in state and control variables subject to mixed control-state constraints}. Optimal Control Applications and Methods 30 (2009).

\bibitem{pontr} L.S. Pontryagin. \textit{Mathematical Theory of Optimal Processes}. John Wiley \& Sons, 1962.

\bibitem{stoer} J. Stoer and R. Bulirsch. \textit{Introduction to Numerical Analysis}. Springer-Verlag New York, 1983.

\bibitem{theo1} A. Rustichini. \textit{Functional differential equations of mixed type: The linear autonomous case}. J. Dynam. Differential Equations, 1 (1989).

\bibitem{theo2} J. Mallet-Paret. \textit{The Fredholm alternative for functional differential equations of mixed type}. J. Dynam. Differential Equations, 11 (1999).

\bibitem{lunel} N. J. Ford, P. M. Lumb, P. M. Lima and M. F. Teodoro. \textit{The numerical solution of forward–backward differential equations: Decomposition and related issues}. Journal of Computational and Applied Mathematics, 234(9), 2745-2756 (2010).

\bibitem{allgower} E. L. Allgower, K. Georg. \textit{Introduction to Numerical Continuation Methods}. Springer Series in Computational Mathematics, Springer. Berlin, Paris. 1990.

\bibitem{Emm} T. Haberkorn, E. Tr\'{e}lat. \textit{Convergence results for smooth regularizations of hybrid nonlinear optimal control problems}. SIAM Journal on Control and Optimization, Society for Industrial and Applied Mathematics, 2011, 49 (4), pp.1498-1522. hal-00519458.

\bibitem{EmmSilva} C. Silva, E. Tr\'{e}lat, \textit{Smooth Regularization of Bang-Bang Optimal Control Problems}. IEEE Transactions On Automatic Control, Vol. 55, No. 11, November 2010.

\bibitem{artEmm} E. Tr\'{e}lat. \textit{Some Properties of the Value Function and its Level Sets for Affine Control System with Quadratic Cost}. Journal of Dynamical and Control Systems (2000).
\end{thebibliography}
\end{document}